\documentclass[runningheads, letter, 12pt]{llncs}
\usepackage{hyperref}

\usepackage{algorithmicx}
\usepackage{algorithm}
\usepackage{algpseudocode}
\usepackage{graphicx}
\usepackage{extarrows}
\usepackage{amssymb,amscd}
\usepackage[all]{xy}
\usepackage{longtable}
\usepackage[latin1]{inputenc}
\usepackage{tikz,arydshln}
\usetikzlibrary{shapes,arrows}

\newcommand{\keywords}[1]{\par\addvspace\baselineskip
\noindent\keywordname\enspace\ignorespaces#1}

\newcommand{\rmv}[1]{}

\newcommand{\Z}{{\mathbb Z}}

\newcommand{\F}{{\mathbb F}}
\newcommand{\C}{{\mathbb C}}
\newcommand{\Q}{{\mathbb Q}}

\newcommand{\fp}{\mathfrak{p}}

\newcommand{\fa}{\mathfrak{a}}

\renewcommand{\O}{{\mathcal O}}
\DeclareMathOperator{\Gal}{\mathrm{Gal}}

\begin{document}
\mainmatter
\title
{Solving the Decision Principal Ideal Problem with Pre-processing}

\author{Jincheng Zhuang \inst{1,2}%
        \and Qi Cheng\inst{3} %
 }
 \institute{
 School of Cyber Science and Technology, Shandong University\\
 Qingdao 266237, China \\
   \and
 State Key Laboratory of Cryptography and Digital Economy Security \\
  Shandong University, Qingdao, 266237, China \\
 Email: {\tt jzhuang@sdu.edu.cn}
  \and 
 School of Computer Science, University of Oklahoma\\
 Norman, OK 73019, USA.\\
 Email: {\tt qcheng@ou.edu}
 }
\maketitle

\begin{abstract}
  The principal ideal problem constitutes a fundamental problem in algebraic
  number theory and has attracted significant attention due to its applications in ideal lattice based cryptosystems. Efficient quantum algorithm has been found to address this problem. The situation is
  different in the classical computational setting. In this work, we delve into
  the relationship between the principal ideal problem and the class field
  computation. We show that the decision version of the problem can be solved
  efficiently if the class group is smooth, after pre-computation has been completed
  to collect information about the Hilbert class field.

\end{abstract}

\section{Introduction}
Given an algebraic number field $K$, let $\O_K$ be the ring of integers of $K$
and $\fa\subset \O_K$ be an ideal. The principal ideal problem (PIP for short) has the following versions:
\begin{enumerate}
  \item (Decision, D-PIP for short) Decide whether $\fa$ is principal.
  \item (Search, S-PIP for short) If $\fa$ is principal, find a generator $g$ such that $\fa=(g)$.
  \item (Optimization, O-PIP for short) If $\fa$ is principal, find a
    short generator $g'$ such that $\fa=(g')$.
\end{enumerate}

The class of principal ideals form the identity element in the class group.
It should be the first question we ask about a given ideal
if the ring of integers of a field is not a principal ideal domain.
It is surprising that all the three versions of PIP are
hard when the degree of a number field is high.
Not only these problems are  fundamental in algebraic
number theory, but they also have relations with other important problems
such as class group discrete logarithm,
and have been studied in cryptography.
For instance, based on the work of Campbell, Groves, and Shepherd~\cite{CGS14},
 Cramer, Ducas, Peikert, and Regev~\cite{CDPR16} showed a reduction
from O-PIP to S-PIP in prime-power cyclotomic fields.
Given a cyclotomic number field of prime power conductor with degree $n$,
Cramer, Ducas and Wesolowski~\cite{CDW17} established a reduction from the
$\gamma-$Ideal Shortest Vector Problem, where the approximate
factor $\gamma=e^{\tilde{O}(\sqrt{n})}$, to S-PIP.
Moreover, the PIP problem has applications in a few
cryptosystems based on ideal lattices such as the public key encryption scheme SOLILOQUY~\cite{CGS14}, the fully homomorphic encryption scheme proposed by Smart and Vercauteren~\cite{SV10} and the multilinear maps proposed by Garg, Gentry and Halevi~\cite{GGH13}.

In the quantum computation model, Hallgren~\cite{Hallgren05} designed a polynomial time algorithm solving S-PIP over constant degree number fields.
Campbel, Groves and Shepherd~\cite{CGS14} described a quantum polynomial time algorithm
in totally real cyclotomic fields.
Biasse and Song~\cite{BS16} proposed a polynomial time algorithm solving S-PIP over arbitrary classes of number fields.

In the classical computation model, a serial of algorithms have been proposed.
A general strategy for solving D-PIP and S-PIP is to utilize the computation
of the full class group, which results in subexponential algorithm~\cite{Fieker01,Biasse14,BF14}.
For special cases, efficient tactics are applied which results in much more
efficient algorithms. For example, some algorithms explore the subfield structure
such as~\cite{BBVLV17,BEFGK17}.

\subsection{Our results and technical overview}

The input size of a PIP problem can be measured by several parameters.
For the field \( K \),
we have  two: the degree of \( K \) denoted by \(
d = [K:\Q] \) and
the discriminant of \( K \) denoted by \( \Delta(K) \).
The field \( K \) can be described by \( O( d \log(\Delta(K)) ) \) bits
of data.
For an input ideal $I$, we can use its log-norm $\log ( N(I) )$
to estimate its bit length.
A PIP algorithm is efficient if it runs
in time \( poly(d, \log \Delta(K), \log N(I)) \).
Our main result can be summarized as follows:

\begin{proposition}
  Suppose that the class group \( Cl(K) \)
  of a number field K is isomorphic to
  \[   \Z/q_1 \times \Z/q_2 \times \cdots \Z/q_t, \]
  where \( q_i \) are prime powers (may not be distinct).
  Let \( s=max_i \{ q_i \}\).
After pre-computation produces \( poly(d, \log(\Delta(K)), s) \) bits
of data, there is an algorithm which decides whether an
ideal $I$  is principal or not,
in time \( poly(d, \log(\Delta(K)), s, \log(N(I))) \),
under a reasonable number theoretical assumption.
\end{proposition}

The main idea to prove the proposition is to use the Hilbert class (sub)field
if the input ideal is prime. For a general ideal,
we rely on a randomized ``ideal switching'' method to
reduce it to the prime ideal case. The switching algorithm
is efficient, assuming a plausible number-theoretic
conjecture 
that was proved in the average case in \cite{deBoer22}.

Algorithms with pre-computation belong to the class of
non-uniform algorithms in theoretical computer science.
They may take subexponential or even exponential time
to produce advice based on part of the input,
and rely on the advice, which should not be very long,
to speed up the computation. In cryptography, it is desirable
for the sake of forward secrecy to
choose cryptosystems without pre-computation  attacks,
since some part of the keys  is fixed and public.
Take the discrete logarithm over finite fields as an example.
In the problem, we are given a prime \( p \), and an element
\( b \) which generates (a large subgroup of) \( ( \F_p )^* \).
For an input \( a \in ( \F_p )^* \), we need to find \( x \) so that
\( a = b^x \pmod{p} \). There is no known exponential-time pre-computation
that can be done on \( p \) and \( b \) to produce advice  so
that the discrete logarithm can be solved efficiently
(although for practical consideration, one can
find an interesting discussion in \cite{BL12}).
Compare it with RSA that can easily be broken by pre-computation
(that finds the factors of the RSA modulus as the advice).
It explains why the Diffie-Hellman
key exchange is more preferable in many situations.

Previous algorithms solving D-PIP, even with
pre-computation,  must run many iterations of
BKZ lattice reductions to achieve a smoothness condition.
Each BKZ reduction  by itself will also take
 sub-exponential time using classical computers.
Furthermore all previous algorithms to solve D-PIP  are not sensitive
to the size of  the class group. In the extreme case, even
the class group is very small (a special case of a smooth group),
previous algorithms will still run in classical
sub-exponential time, whereas our new non-uniform algorithm has complexity
depending on the size of class group:

\begin{corollary}
  Let \( I \) be an ideal in a number field \( K \) with norm \( N(I) \).
  Assuming Conjecture~\ref{chanceofprime},
  there is a non-uniform algorithm, taking
  \( poly(d, \log(\Delta(K)), |Cl(O_K)|) \) bits of data as
  advice, running in time
  \[ poly(d, \log(\Delta(K)), |Cl(O_K)|, \log(N(I))), \]
  and deciding whether the input ideal \( I \) is principal or not.
\end{corollary}

The main drawback of our algorithm is that it does not
calculate generators of ideals after they are found
to be principal. In the other words, it does not solve the S-PIP.

\section{Preliminaries}
\subsection{Number fields}
A number field $K$ is a finite extension of $\Q$ with degree denoted by
$d=[K:\Q]$. Every number field $K$ is a simple extension of $\Q$ such that
$K=\Q(\theta)$ for certain $\theta\in\C$. An element $\alpha\in K$ is called
integral if it satisfies a monic polynomial in $\Z[x]$.
The set of integral elements of $K$ forms a ring, denoted as $\O_K$.
This ring of integers  is Dedekind, thus it enjoys the unique
factorization of ideals. Let $\fa$ be an ideal of $\O_K$, then
we have
\[
\fa=\fp_1^{e_1}\cdots \fp_g^{e_g},
\]
where $\fp_i$ is a prime ideal, $e_i\in \Z$ are non-negative and the factorization is unique up to ordering.

A fractional ideal $I$ of $\O_K$ is a nonzero $\O_K$-submodule
of $K$ such that $\alpha I\subset \O_K$ for certain $\alpha\in K^*$. The fractional ideal also enjoys the property of unique
factorization. That is, given a fractional ideal $I$, then we have a unique factorization up to ordering
\[
I=\fp_1^{e_1}\cdots \fp_s^{e_s},
\]
where $\fp_i$ is a prime ideal of $\O_K$, $e_i\in \Z$.

Let $I_K$ be the group of  fractional ideals
under the ideal multiplication.
Let $P_K=\{\alpha\O_K:\alpha\in K^*\}$ be the subgroup of principal ideals. The class group of $\O_K$ is defined as
\[
Cl(\O_K)=I_K/P_K,
\]
the cardinality of which measures how far $\O_K$ is away from
unique factorization domain. We refer the reader to~\cite{Marcus18}\cite{SD12} and references therein for further properties of number fields.

\subsection{Hilbert class field}

Let $K$ be a number field.
The Hilbert class field $H$ of $K$ is defined as its maximal
abelian unramified extension. If \( H \) is known,
one can quickly tell whether a prime ideal in \( K \)
is principal or not.

\begin{theorem}[\cite{SD12}, Section 17]
  Let $H$ be the Hilbert class field of a number field $K$,
  and let $\fp$ be a prime ideal of $K$. Then $\fp$ is principal in \( K \)
  if and only if $\fp\O_H$ splits completely in \( H \).
\end{theorem}

Hilbert class fields encode information about the class group, that is,
\[ Gal(H/K) = Cl(K).  \]
If the class group is large, then the Hilbert class field has high degree.
It may be sufficient to describe all its low degree subfields
if those subfields generate the Hilbert class field.

Fieker~\cite{Fieker01} described an algorithm to compute
the Hilbert class field, and presented codes and numerical examples.
He did not analyze the time complexity of his algorithm, but it
is clear from the paper that the most expensive step in the algorithm
is to compute the class group, and the rest of his algorithm has
complexity depending on the degree of the Hilbert class field.
Eisentr{\"{a}}ger and Hallgren~\cite{EH10} gave an analysis to
the algorithm of constructing a small degree subfield of a Hilbert class field,
which is a simplified version of the Algorithm~5.2.14 described in~\cite{Coh00}.
They reduce the Hilbert class field computation problem to various problems in algorithmic number
theory. We will employ these algorithms of computing Hilbert class field in our pre-computation.

\subsection{Main ideas in the previous results}

For quadratic imaginary fields, deciding the PIP can be solved
efficiently by reducing quadratic forms~\cite{Cox89}.
When the degree of fields become high,
the PIP can be solved in sub-exponential
time by lattice reduction, or quantum polynomial time. Here we
review the ideas in the classical sub-exponential time algorithm.
First it is known that under the Generalized Riemann Hypothesis
(GRH), a class group can be generated by a small number
of ideals. Let $C$ be a sub-exponential function
in $d$ and \( \log \Delta(K) \). Let
\[ \fp_1, \fp_2, \cdots, \fp_C \]
be the prime ideals with small norms, which by the GRH, generate
the class group.

A vector \( (e_1, e_2, \cdots, e_C) \)
is called a relation if \( \prod_{1\leq i\leq C} \fp_i^{e_i} \)
is principal.
All the relations form a \( \Z \)-modulo which will be denoted
by P. Once we collect enough number of relations,
we find a base for P and complete an exact sequence:
\[ 0 \rightarrow P \rightarrow \Z^C \rightarrow Cl(K) \rightarrow 0.  \]
It allows us to find the generators for \( Cl(K) \)
by computing the Smith Normal Form of the relation matrix.

If the input ideal \( I \) has large norm,
we need to reduce the size of an ideal modulo principal ideals.
In other words, we find an ideal which is in the same (or the inverse)
class  of \( I \), and whose norm depends only on \( d \) and
\( \log \Delta(K) \), and is independent to \( I \).
This step is called the close principal multiple problem (CPM for short) in \cite{CDW17}.
The algorithm suggested in many papers such as
\cite{BEFGK17} is to embed the ideal into a lattice and
to find a short element \( v \) in $I$ by lattice reduction.
Then the algorithm proceeds to compute the ideal \( I' \)
so that \( (v)=I' I \). We will run the BKZ reduction with block
size \( \beta \) so \( v \) has
high quality, which guarantees that the determinant of \( I' \)
is independent to the norm of \( I \).
More precisely, if the block size is   \( \beta \)
 we have
\[ N(I') \leq \beta^{d^2/\beta} \sqrt{|\Delta(K)|}. \]
Note that running LLL ( \( \beta=2 \) )
may work if the degree \( d \) is small \cite{CDO97}, but
is not good enough if the degree of the field is
high. We need to run BKZ with large
block size around \( \beta = \sqrt{d} \) so that
\( L_{1/2}(N(I')) \) remains sub-exponential in
\( d \) and \( \log (\Delta(K)) \).
Here  \[ L_{1/2}(x) = exp(\tilde{O}( (\log x)^{1/2} )). \]
One hopes that \( I' \)
 is smooth, namely, it can be written
as a product of prime ideals from \( \fp_1, \fp_2, \cdots, \fp_C \).
The probability is about 1 over a sub-exponential function
in \( d \) (and \( \log \Delta(K) \)),
assuming a few number theory conjectures.
One can see that all the previous algorithms deciding
whether an ideal is principal or not
requires many expensive lattice reductions,
each of which runs in a sub-exponential time.
For details see~\cite{BEFGK17}.
In this paper we will apply new ideas
to decrease the complexity of the D-PIP.
We note that the quantum polynomial time algorithm can solve
both D-PIP and S-PIP.

\section{Non-Uniform Algorithm for D-PIP}
In this paper, we introduce a parameter  to measure
the smoothness of a finite abelian group. It is a generalization
of the smoothness of integers.
A finite abelian group \( G \) can be decomposed as
\[ G = \Z/q_1 \times \Z/q_2 \times \cdots \times \Z/q_t,  \]
where \( q_i \)'s are prime powers (some of them may be equal).
We say the group is \( s \)-smooth if all the \( q_i \)'s are
less than \( s \).
We consider the situation when $\Gal(H/K)$ is smooth,
and reduce the decision PIP to the polynomial factorization
over finite fields.

The relation of the prime ideal splitting in the composite number field
and its subfields is described in the following theorem.
\begin{theorem}[\cite{Marcus18}, Theorem 31]
  Let $E_1,\cdots,E_t$ be Galois extensions of a number field $K$,
  $L$ be the composite field of $E_1,\cdots,E_t$, $\fp$ be a prime ideal
  of $\O_K$. Then
  \[
  \fp\O_L \ \text{splits completely} \iff \fp\O_{E_i} \text{splits completely for} \  1\leq i\leq t.
  \]
\end{theorem}

As an application, we have the following corollary.

\begin{corollary}
Let $H$ be the Hilbert class field of a number field $K$, such that
\[ \Gal(H/K)=\Z/q_1 \times \Z/q_2 \times \cdots \times \Z/q_t, \]
where $q_i$ are prime powers (may not be distinct).
Assume \[ H=E_1E_2\cdots E_t \] such that $|\Gal(E_i/K)|=q_i$.
Let $\fp$ be a prime ideal of $\O_K$. Then
  \[
  \fp \ \text{is principal} \iff \fp\O_{E_i} \text{splits completely for} \  1\leq i\leq t.
  \]
\end{corollary}

\begin{figure}[h]
\[
\xymatrix{
     &     & H \ar@{-}[lld] \ar@{-}[ld]\ar@{-}[rd]\ar@{-}[rrd]  & & \\
 E_1 \ar@{-}[rrd]|-{q_1} & E_2\ar@{-}[rd]|-{q_2} &  & \cdots\ar@{-}[ld] & E_t\ar@{-}[lld]|-{q_t} \\
    &     &  K  & & \\
}
\]
\caption{Extension diagram in the smooth case}\label{Fig1}
\end{figure}
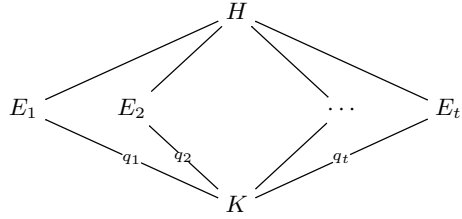

Pictorially, we have the diagram as shown in Fig.~\ref{Fig1}.


\subsection{Pre-computation and algorithm for prime ideals}

In this paper, we denote the discriminant of a polynomial \( f \)
by \( \Delta(f) \). For a number field \( L \) that is an extension
of \( K \), we denote the relative discriminant of \( L \)
with respect to \( K \) by \( \Delta_K (L) \). If \( K = \Q \),
we will omit the subscription sometimes.

The pre-computation stage will compute the fields \( E_1, E_2, \cdots, E_t.  \)
More precisely we will find all the monic polynomials
\( f_1, f_2, \cdots, f_t \) in \( O_K [x] \)
 so that for all \( 1\leq i\leq t \)
\[E_i = K[x]/(f_i(x)). \]
Such algorithm is described in \cite{Fieker01}.
In general, a prime ideal \( \fp \subseteq O_K\) dividing \( p\O_K \)
is principal iff for all \( 1\leq i\leq t \), \( f_i \) splits completely in
\( O_K [x]/(\fp) \), which is a finite field with characteristic \( p \).
Special care needs to be taken for prime factors of discriminant of \( f_i \).
For those primes, we decide whether they are principal in the pre-computation
stage, and if they are principal, we put them in a set denoted by \( S \).

\begin{proposition}\label{advice}
  The pre-computation will produce the polynomials
  \( f_1, f_2, \cdots, f_t \) and a set of
  prime ideals
  \[ S = \{ \fp | \fp {\ is\ principal,\ and\ } \fp {\ divides\ }
    \Delta(f_i) {\ for\ some\ }1\leq i\leq t\}  \]
  in quantum polynomial time, or classic sub-exponential time.
  And the size of advice (including \( f_i \) and \( S \))
  is \( poly(d, \log (\Delta(K)), \max_i \{ q_i\}) \) bits.
\end{proposition}

To bound the size of the advice in Proposition~\ref{advice}, we
first observe that the relative discriminant
\( \Delta_K(E_i) \) of \( E_i \)
over \( K \) is 1 since \( E_i \) is unramified over \( K \).
So the absolute discriminant
of \( E_i \) is
\[ \Delta ( E_i) = \Delta(K)^{q_i},  \]
according to the relative discriminant formula described
for instance in~\cite{FT93}(Proposition III.2.15).
Hence we only need \( poly(q_i, \log \Delta(K)) \)
space to describe the field \( E_i \).
Secondly if we examine the algorithms
to construct the field \( E_i \),
they involve finding ideals \( I \) with order \( q_i \)
in the class group,
and finding an element \( u \) so that \( I^{q_i} = (u)\).
The polynomial \( f_i \) can be derived from
 the minimal polynomial of \( ( \epsilon u )^{1/q_i} \) over K,
where \( \epsilon \) is a unit in \( E_i \).
The number of  prime ideal factors of \( \Delta(f_i) \)
is bounded from above by \( \log |\Delta(f_i)| \).
And finding them
takes quantum polynomial time  or classic sub-exponential time.
All the data may require costly computation, but it can
be written down succinctly.

Denote the input ideal by \( I \).
We describe the algorithm for deciding whether a given prime ideal is principal
in Algorithm~\ref{Prime-Ideal-D-PIP}.

\begin{algorithm}[htb]
  \caption{Prime-Ideal-D-PIP: deciding whether a given prime ideal is principal}
  \label{Prime-Ideal-D-PIP}
  \begin{algorithmic}[1]
    \Require
      Advice: \( f_1, f_2, \cdots, f_t \) and \( S \),
      and a prime ideal \( I \)
    \Ensure
      ``Yes'' if \( I  \) is principal, ``No'' otherwise
     \For { \( 1\leq i \leq t \) }
        \If {\( I | \Delta(f_i): \)}
          \If {\( I \in S \)}
              \State \Return ``Yes''
           \Else
              \State \Return ``No''
          \EndIf
         \EndIf
     \EndFor
     \For {\( 1\leq i \leq t \)}
        \If { \( f_i \pmod{I} \) does not split completely }
           \State \Return ``No''
        \EndIf
     \EndFor
     \State \Return ``Yes''
  \end{algorithmic}
\end{algorithm}

\subsection{Algorithm for general ideals}

For a general ideal that may not be prime, we  reduce
the D-PIP problem to the prime ideal case. Let \( I \)
be the input ideal. Let \( b_1, b_2, \cdots, b_d \)
be its $\Z$-base (assumed to be LLL-reduced).
Let $B$ be a positive integer.
We find \( d \) integers \( r_1, r_2, \cdots, r_d \)
uniformly random in the range \( [-B,B] \).
Let \( r = r_1 b_1 + r_2 b_2 + \cdots + r_d b_d \).
Compute the ideal \( I' \) so that
\( (r) = I' I \).
What is the chance that \( I' \) is a prime ideal?
Heuristically \( I' \) behaves like a random ideal
with bounded determinants
in the class of \( [I]^{-1} \) when \( B \)
is large enough, so the chance is high. Here we use \( [I] \)
to denote the ideal class containing \( I \).

\begin{conjecture}\label{chanceofprime}
  Let \( B = 2^d \Delta(K)  \). Denote
  \[
  N = \left| \left\{ (r_1,r_2, \cdots, r_d)\in \Z^d:
      \begin{matrix}
    \forall 1\leq i\leq d, -B\leq r_i\leq B, \\
    {\rm and\ } ( \sum_{i=1}^d r_i b_i )/I \\
        {\rm \ is\ a\ prime\ ideal}
        \end{matrix}
        \right\} \right|.
  \]

  We have
  \[ \frac{N}{(2B)^d}
    > \frac{1}{poly( d,\log(B))}.\]
\end{conjecture}

In essence, the conjecture asserts that the probability that \( I' \)
is a prime ideal should be non-negligible when \( B \)
is large enough.
According to the Landau prime ideal theorem, the number of prime ideals
with norm bounded by $T$, denoted as $\pi_K(T)$, is approximately
\[
\pi_K(T)\sim \frac{T}{\log T}.
\]
In Chapter 6 of the Phd dissertation \cite{deBoer22},
the conjecture is proved if the ideal \( I \)
is chosen uniformly in random from the Arakelov class group.
We also provide some numerical evidence in Section~\ref{numevidence}.
The computation data shows that it does not take long to
find prime ideals even for small \( B \).

If the conjecture is true, then we can find a prime ideal
\( I' \) quickly. And \( I' \)
is principal iff \( I \) is principal. We then run
the algorithm of the prime case.
We describe the algorithm for deciding whether a given general ideal is principal
in Algorithm~\ref{General-Ideal-D-PIP}.
This is a randomized algorithm which should solve the problem with
non-negligible probability, assuming Conjecture~\ref{chanceofprime}.

\begin{algorithm}[htb]
  \caption{General-Ideal-D-PIP: deciding whether a given general ideal is principal}
  \label{General-Ideal-D-PIP}
  \begin{algorithmic}[1]
    \Require
      an ideal \( I \) given by its \( \Z \)-base \( b_1, b_2, \cdots, b_d \)
    \Ensure
      ``Yes'' if \( I  \) is principal, ``No'' otherwise
    \Repeat \For { \( 1\leq i \leq d \) }
               \State \( r_i \overset{ \text{random} }{\longleftarrow} [-B,B] \)
            \EndFor
           \State \( r = \sum_{1\leq i\leq d} r_i b_i \)
           \State \( I' = (r)/I \)
     \Until { \( I' \) is a prime ideal }
     \State \Return Prime-Ideal-D-PIP($I'$)
  \end{algorithmic}
\end{algorithm}

\section{Numerical Examples}\label{numevidence}
We compute some examples using Pari~\cite{PARI} and SageMath~\cite{Sagemath} to illustrate the algorithm.

\begin{example}
First we present a toy example. Let \( K=\Q(\sqrt{-5}) \).
The Hilbert class field is \( H = K[\sqrt{-1}] = K[x]/(x^2+1) \).
Here the polynomial \( x^2+1 \) is taken as ``advice''
in our algorithm.
Its discriminant  is \( -4 \), whose prime factor \( 2 \)
ramifies in \( K \) i.e. \( (2) = (2, 1+\sqrt{-5})^2 \).
The ideal \(  (2, 1+\sqrt{-5})\) is not principal,
and it should not be put into the set \( S \) in Proposition~\ref{advice}.
One can verify that for a prime ideal \( \fp \)
in \( K \) other than \( (2, 1+\sqrt{-5}) \),
\( \fp \) is principal if and only if \( x^2+1  \)
splits into two distinct linear factors \( \pmod{\fp} \).
For example, if a rational prime $p$ is inert in $K$,
then  $(p)$ is principal, $x^2+1$  will of course
split in \( \F_{p^2} = O_K/(p) \).
It can also be verified that \( H = K[\sqrt{5}] = K[x]/(x^2-x-1).  \)
The algorithm will give the same answer no matter
we use \( x^2+1 \) or \( x^2-x-1 \).
To see that, we only need to examine  rational prime \( p \)
that splits into a product of two prime ideals \( \fp_1 \fp_2 \) in \( K \).
The two polynomials $x^2 + 1$ and $x^2-x-1$ will have the same factorization
pattern \( \pmod{\fp_1} \), since their roots are related due to
\( \sqrt{-5} = \sqrt{-1}\sqrt{5} \),
and \( \sqrt{-5} \) is congruent to an integer \( \pmod{\fp_1} \).
\end{example}

\begin{example}
Let \( K=\Q(\zeta_{180})=\Q[z]/(z^{48}+z^{42}-z^{30}-z^{24}-z^{18}+z^6+1)\), whose class number is
$75$ and class group $Cl(K)\cong \Z_3\times\Z_5\times\Z_5$.
The Hilbert class field is the composition of \(E_1=K[x]/(f_1(x)), E_2=K[x]/(f_2(x)), E_3=K[x]/(f_3(x))\), where
\[
\begin{aligned}
f_1(x) = &\ x^3 + (-z^{45} + 2z^{42} - z^{39} + z^{36} - z^{30} + z^{27} - z^{24} - z^{18} \\
         &\ + z^{15} - z^{12} + 1), \\
f_2(x) = &\ x^5 + (-5z^{44} - z^{43} - 2z^{40} - z^{37} - 2z^{35} + 5z^{32}  + 5z^{26} \\
         &\ + 6z^{25} - z^{20} + z^{19} + z^{15} + z^{13} - 4z^{10} - 5z^8 + 6z^5\\
         &\   - 5z^2 - z - 1),\\
f_3(x) = &\ x^5 + (4z^{45} + 47z^{44} + 28z^{43} + 44z^{40} + 28z^{37} - 34z^{35}  \\
         &\ - 47z^{32}- 21z^{30} - 47z^{26} + 21z^{25}- 13z^{20} - 28z^{19} \\
         &\   - 44z^{15} - 28z^{13}- 4z^{10} + 47z^8 + 47z^5 + 47z^2 \\
         &\  + 28z - 28). \\
\end{aligned}
\]
So to decide whether a prime ideal is principal or not, one
only needs to factor three polynomials of degree 3, 5 and 5
respectively, rather than reducing many lattices of dimension 48.
\end{example}

\begin{example}
Let \( K=\Q(\zeta_{64})=\Q[z]/(z^{32}+1)\).
Consider the ideal \(I=(187,34z^{16} - 85z^8 - 33z^4 + 54)\), which is not prime and has reduced basis
as shown in Table~\ref{Table-Basis}.

\begin{table}
\caption{Reduced basis of ideal I}\label{Table-Basis}
\centering
\resizebox{0.9\columnwidth}{!}{
\begin{tabular}{|c|c|c|c|c|c|c|c|}
  \hline
  Basis & Element & Basis & Element \\
  \hline
$b_1$ & $z^{24}+z^{20}+z^{16}+3z^{12}+4z^8-z^4+2$   & $b_2$   & $-z^{28}-z^{24}-z^{20}-3z^{16}-4z^{12}+z^8-2z^4$ \\
$b_3$   & $-z^{28}-z^{24}-3z^{20}-4z^{16}+z^{12}-2z^8+1$ &$b_4$   & $3z^{24}-3z^{20}-z^{16}+2z^{12}+3z^4-1$  \\
$b_5$   & $-z^{28}-3z^{20}+3z^{16}+z^{12}-2z^8-3$   & $b_6$   &$-3z^{28}+z^{24}+3z^{16}-3z^{12}-z^8+2z^4$ \\
$b_7$   & $z^{28}-2z^{24}+z^{16}+z^{12}+z^8+3z^4+4$   & $b_8$   & $2z^{28}-z^{20}-z^{16}-z^{12}-3z^8-4z^4+1$ \\
$b_9$ & $z^{25}+z^{21}+z^{17}+3z^{13}+4z^9-z^5+2z$ &$b_{10}$ & $-z^{29}-z^{25}-z^{21}-3z^{17}-4z^{13}+z^9-2z^5$ \\
$b_{11}$   & $-z^{29}-z^{25}-3z^{21}-4z^{17}+z^{13}-2z^9+z$&$b_{12}$ & $3z^{25}-3z^{21} -z^{17}+2z^{13}+3z^5-z$  \\
$b_{13}$   & $-z^{29}-3z^{21}+3z^{17}+z^{13}-2z^9-3z$   & $b_{14}$   &$-3z^{29}+z^{25}+3z^{17}-3z^{13}-z^9+2z^5$ \\
$b_{15}$   & $z^{29}-2z^{25}+z^{17}+z^{13}+z^9+3z^5+4z$ & $b_{16}$ & $2z^{29}-z^{21}-z^{17}-z^{13}-3z^9-4z^5+z$ \\
$b_{17}$ & $z^{26}+z^{22}+z^{18}+3z^{14}+4z^{10}-z^6+2z^2$&$b_{18}$& $-z^{30}-z^{26}-z^{22}-3z^{18}-4z^{14}+z^{10}-2z^6$ \\
$b_{19}$   & $-z^{30}-z^{26}-3z^{22}-4z^{18}+z^{14}-2z^{10}+z^2$ &$b_{20}$   & $3z^{26}-3z^{22}-z^{18}+2z^{14}+3z^6-z^2$  \\
$b_{21}$   & $-z^{30}-3z^{22}+3z^{18}+z^{14}-2z^{10}-3z^2$   & $b_{22}$   &$-3z^{30}+z^{26}+3z^{18}-3z^{14}-z^{10}+2z^6$ \\
$b_{23}$   & $z^{30}-2z^{26}+z^{18}+z^{14}+z^{10}+3z^6+4z^2$   & $b_{24}$   & $2z^{30}-z^{22}-z^{18}-z^{14}-3z^{10}-4z^6+z^2$ \\
$b_{25}$ & $z^{27}+z^{23}+z^{19}+3z^{15}+4z^{11}-z^7+2z^3$ &$b_{26}$ & $-z^{31}-z^{27}-z^{23}-3z^{19}-4z^{15}+z^{11}-2z^7$ \\
$b_{27}$   & $-z^{31}-z^{27}-3z^{23}-4z^{19}+z^{15}-2z^{11}+z^3$ &$b_{28}$ & $3z^{27}-3z^{23}-z^{19}+2z^{15}+3z^7-z^3$  \\
$b_{29}$   & $-z^{31}-3z^{23}+3z^{19}+z^{15}-2z^{11}-3z^3$   & $b_{30}$   &$-3z^{31}+z^{27}+3z^{19}-3z^{15}-z^{11}+2z^7$ \\
$b_{31}$   & $z^{31}-2z^{27}+z^{19}+z^{15}+z^{11}+3z^7+4z^3$ & $b_{32}$ &
$2z^{31}-z^{23}-z^{19}-z^{15}-3z^{11}-4z^7+z^3$ \\
  \hline
\end{tabular}
}
\end{table}

Testing the technique of random ideal switching in Algorithm~\ref{General-Ideal-D-PIP},
the average number of tests needed before obtaining a prime ideal $I'$
 is summarized in the Table~\ref{Ideal-Switch}.
\begin{table}
\caption{Average number of ideal switches until prime}\label{Ideal-Switch}
\centering
\begin{tabular}{|c|c|}
  \hline
  Bound $B$ &  Average number of  ideal switches until prime\\
  \hline
  5 &   20   \\
  10 &   26  \\
  20 &  32  \\
  \hline
\end{tabular}
\end{table}
\end{example}

\section{Conclusion and Open Problems}
We solve the decision version of PIP efficiently when the class group
is smooth, assuming that pre-computation has been done to
collect information about the Hilbert class field.
Our results indicate that there is a significant gap in
computation complexity between the decision version and
the search version of the PIP.
The main open problem is to close the gap, i.e. to
design an efficient algorithm that can find
a generator of a principal ideal when the class group is smooth.

Another interesting open problem is to calculate
the low degree subfield of the Hilbert class field,
when the Hilbert class field has other sub-fields of high degrees.
For example, the number field \( F = \Q(\zeta_{512}) \)
has a large class group, whose order is divisible by 17, since
the relative class number
\[  h^{-}=17\times 21121\times 76 532353\times 29102880226241\times p_{28}, \]
where \( p_{28} \) is prime number of 28 decimal digits.
So it would take too much space to describe the Hilbert class field
\( H \) of \( F \), but is it possible to find a
monic irreducible polynomial \( f_1 \in O_F[x] \) so
that \( F[x]/ (f_1) \) is the sub-field of \( H \)
of degree \( 17 \) over \( F \) in reasonable amount of time?

\pagestyle{plain}
\bibliographystyle{splncs04.bst}
\bibliography{PIP-CFT}
\end{document}